\documentclass[11pt]{amsart}

\usepackage{pstricks,graphicx,pst-node,amsmath,amssymb,mathdots}
\usepackage{hyperref}
\newtheorem{theorem}{Theorem}
\newtheorem{lemma}[theorem]{Lemma}

\newtheorem{corollary}[theorem]{Corollary}

\theoremstyle{definition}
\newtheorem{definition}[theorem]{Definition}

\theoremstyle{remark}

\begin{document}


\title[A bijection between permutations and DPPs with no special parts]{A direct bijection between descending plane partitions with no special parts and permutation matrices}


\author[J. Striker]{Jessica Striker} %
\address{ }
\email{jessica@math.umn.edu}


\date{\today}



\begin{abstract}
We present a direct bijection between descending plane partitions with no special parts and permutation matrices. This bijection has the desirable property that the number of parts of the descending plane partition corresponds to the inversion number of the permutation. 
Additionally, the number of maximum parts in the descending plane partition corresponds to the position of the one in the last column of the permutation matrix. We also discuss the possible extension of this approach to finding a bijection between descending plane partitions and alternating sign matrices.
\end{abstract}

\maketitle


\section{Introduction}
\label{s:introduction}

\thispagestyle{empty}

Alternating sign matrices have been objects of much inquiry over the past several decades with the most recent exciting development being the proof of the Razumov-Stroganov conjecture~\cite{RAZ_STROG}. Descending plane partitions are objects equinumerous with alternating sign matrices, but there is no bijective proof known. In this paper, we find a direct bijection between these two sets of objects in the simplest special case.

\begin{definition}
A \emph{descending plane partition} (DPP) is an 
array of positive integers $\left\{ a_{i,j}\right\}$ with $i\leq j$ (that is, with the $i$th row indented by $i-1$ units) with weak decrease across rows, strict decrease down columns, and in which the number of parts in each row is strictly less than the largest part in that row and is greater than or equal to the largest part in the next row.
\end{definition}

\begin{figure}[hbp]
$\begin{array}{cccccccc}
a_{1,1} & a_{1,2} & a_{1,3} & \hdotsfor{4} & a_{1,\lambda_1}\\
        & a_{2,2} & a_{2,3} & \hdotsfor{3} & a_{2,\lambda_2} &\\
        &         & \ddots & & & \iddots & &\\
        &         &         & a_{\ell,\ell} & \cdots & a_{\ell,\lambda_\ell} & &
\end{array}$
\caption{The general form of a descending plane partition.}
\label{fig:othernew}
\end{figure}


\begin{definition}
A descending plane partition is \emph{of order $n$} if its largest part is at most $n$.
\end{definition}

\begin{definition}
A \emph{special part} of a descending plane partition is a part $a_{i,j}$ such that \linebreak $a_{i,j}\le j-i$.
\end{definition}

\begin{figure}[htbp]
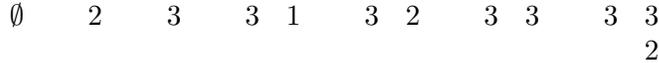

$
\begin{array}{c}
\emptyset\\
~
\end{array}
\hspace{.5cm}
\begin{array}{c}
2\\
~
\end{array}
\hspace{.5cm}
\begin{array}{c}
3\\
~
\end{array}
\hspace{.5cm}
\begin{array}{cc}
3 & 1\\
& 
\end{array}
\hspace{.5cm}
\begin{array}{cc}
3 & 2\\
& 
\end{array}
\hspace{.5cm}
\begin{array}{cc}
3 & 3\\
& 
\end{array}
\hspace{.5cm}
\begin{array}{cc}
3 & 3 \\
&2
\end{array}
$
\caption{The seven descending plane partitions of order 3.}
\label{fig:newname}
\end{figure}

See Figure~\ref{fig:othernew} for the general form of a DPP and Figure~\ref{fig:newname} for the seven DPPs of order $3$. The only DPP in Figure~\ref{fig:newname} with a special part is 3~1. The 1 is a special part since $1=a_{1,2}\leq 2-1$.

Though the definition of descending plane partitions seems a bit contrived, the history behind the counting formula and connections with other combinatorial objects make descending plane partitions interesting objects to study. 

In 1982 Mills, Robbins, and Rumsey proved that DPPs with largest part at most $n$ are counted by the following expression~\cite{MRRANDREWSCONJ}.
\begin{equation}
\label{eq:prod}
\prod_{j=0}^{n-1} \frac{(3j+1)!}{(n+j)!}
\end{equation}
Additionally, they proved that the generating function for DPPs 
of order $n$
weighted by $q$ to the sum of the entries 
is equal to the $q$-ification of (\ref{eq:prod}),
\begin{equation}
\label{eq:qprod}
\prod_{j=0}^{n-1} \frac{(3j+1)!_q}{(n+j)!_q},
\end{equation}
where $k!_q=(1+q)(1+q+q^2)\cdots (1+q+q^2+\cdots +q^{k-1})$~\cite{MRRANDREWSCONJ}. 

Then in 1983 Mills, Robbins, and Rumsey conjectured that $n\times n$ alternating sign matrices were also counted by (\ref{eq:prod})~\cite{MRRASMDPP}.
\begin{definition}
An \emph{alternating sign matrix} (ASM) is a square matrix with entries 0, 1, or $-1$ whose rows and columns each sum to 1 and such that the nonzero entries in each row and column alternate in sign.
\end{definition}

\begin{figure}[htbp]
\scalebox{.85}{
$
\left( 
\begin{array}{rrr}
1 & 0 & 0 \\
0 & 1 & 0\\
0 & 0 & 1
\end{array} \right)
\left( 
\begin{array}{rrr}
0 & 1 & 0 \\
1 & 0 & 0\\
0 & 0 & 1
\end{array} \right)
\left( 
\begin{array}{rrr}
1 & 0 & 0 \\
0 & 0 & 1\\
0 & 1 & 0
\end{array} \right)
\left( 
\begin{array}{rrr}
0 & 1 & 0 \\
1 & -1 & 1\\
0 & 1 & 0
\end{array} \right)
\left( 
\begin{array}{rrr}
0 & 1 & 0 \\
0 & 0 & 1\\
1 & 0 & 0
\end{array} \right)
\left( 
\begin{array}{rrr}
0 & 0 & 1 \\
1 & 0 & 0\\
0 & 1 & 0
\end{array} \right)
\left( 
\begin{array}{rrr}
0 & 0 & 1 \\
0 & 1 & 0\\
1 & 0 & 0
\end{array} \right)
$
}
\caption{The seven $3\times 3$ ASMs.}
\label{fig:3x3asms}
\end{figure}

See Figure~\ref{fig:3x3asms} for the seven $3\times 3$ ASMs. It is clear from the definition that permutation matrices are the alternating sign matrices with no $-1$ entries.
The alternating sign matrix conjecture, as the conjecture of Mills, Robbins, and Rumsey came to be called, stayed open for over a decade and was proved independently by Zeilberger~\cite{ZEILASM} and Kuperberg~\cite{KUP_ASM_CONJ} in 1996.
Many open problems remain concerning DPPs and ASMs, including that of constructing an explicit bijection between them. 
Not even a weight has yet been found on ASMs which yields (\ref{eq:qprod}) as the generating function. 

In Section~\ref{s:sn} of this paper, we present a simple bijection between descending plane partitions of order $n$ with no special parts and $n\times n$ permutation matrices, which are the alternating sign matrices with no $-1$ entries. This bijection has the property that the number of parts of the descending plane partition corresponds to the inversion number of the permutation. 
Also, if there are $k$ parts equal to $n$ in the DPP then there will be a 1 in row $n-k$ of the last column of the corresponding ASM. 

%
These two properties of the bijection make it a bijective solution in the special case $m=0$ of Conjecture 3 in~\cite{MRRASMDPP} which states that the number of DPPs of order $n$ with $p$ parts, $k$ parts equal to $n$, and $m$ special parts  is equal to the number of $n\times n$ ASMs with (generalized) inversion number $p$, $m$ entries equal to $-1$, and a $1$ in column $k+1$ of the first row (or equivalently, a $1$ in row $n-k$ of the last column). 
Thus, the bijection of this paper may be a good first step toward a bijective solution of this conjecture for arbitrary $m$.

Recently, Behrend, Di Francesco, and Zinn-Justin proved this conjecture for all $m$ in a non-bijective way by finding  determinant expressions for each generating function and then proving the equality of these determinants~\cite{ZINNDPP}. Also, in~\cite{LALONDE}, Lalonde gave a bijective proof in the case $m=1$ and mentions the case $m=0$, but does not give the details of the bijection in the case of no special parts.

Additionally, in~\cite{AYYER_ASM_DPP} Ayyer found a different bijection between DPPs of order $n$ with no special parts 
and permutations of $n$, which does not map DPPs with $p$ parts to permutations with $p$ inversions, but rather maps DPPs with $\ell$ rows to permutations with $\ell$ ascents. Already, this bijection seems complicated and it is unclear how one might extend it to all DPPs and ASMs. 

The proof of our bijection is based entirely on the definitions and does not rely on the bijection of~\cite{AYYER_ASM_DPP}.
It is not obvious how to extend our approach to find a bijection between alternating sign matrices with $m$ entries equal to $-1$ and descending plane partitions with $m$ special parts, but in Section~\ref{s:fullbij} we discuss some promising aspects and difficulties encountered. Perhaps the combined insights from both these bijections will provide a way forward.

\section{A bijection on permutations}
\label{s:sn}

In this section we give a bijection between descending plane partitions
of order $n$ with no special parts, $p$ parts, and $k$ parts equal to $n$
and $n\times n$ permutation matrices with $p$ inversions and a $1$ in row $n-k$ of the last column. We will need the following lemma.

\begin{lemma}
\label{thm:dpppart}
There is a natural part-preserving bijection between descending plane partitions of order $n$  with no special parts 
and partitions with largest part at most $n$ and with at most $i-1$ parts equal to $i$ for all $i\leq n$.
\end{lemma}
The proof of Lemma~\ref{thm:dpppart} is slightly technical, but the bijection map is very simple, so we first state the map then prove that it is a well-defined bijection. To map from the DPPs to the partitions, take all the parts of the DPP and put them in one row in decreasing order. To map from the partitions to the DPPs, insert the parts of the partition into the shape of a DPP in decreasing order, putting as many parts in a row as possible and then moving on to the first position in the next row whenever the next part to be inserted would be forced to be special if added to the current row. So this bijection is simply a rearrangement of parts.
\begin{proof}[Proof of Lemma~\ref{thm:dpppart}]
We start with a DPP $\delta$ of order $n$ with no special parts. We then take all the parts of the DPP and put them in one row in decreasing order to form a partition. It is then left to show that $\delta$ has at most $i-1$ parts equal to $i$ for all $i\le n$. Suppose there exists an integer $i$ with $1\le i\le n$ such that there are at least $i$ parts in $\delta$ equal to $i$. Since there is strict decrease down columns in $\delta$ it follows that there can be at most one $i$ in each column. Since $\delta$ has no special parts, $i$ must be greater than its column minus its row. So $i$ can appear no further right than entry $a_{1,i}$. So to have $i$ parts equal to $i$ in $\delta$, there must be an $i$ in each column from 1 to $i$. The tableau is shifted, so this means that there must be an $i$ in entry $a_{1,1}$. Thus, by the fact that the columns have to be strictly decreasing, all the $i$'s must appear in the first row. But the number of parts in each row must be strictly less than the largest part in that row, so we cannot have the first row of $\delta$ consisting of $i$ parts equal to $i$. Thus there are at most $i-1$ parts equal to $i$.

To map in the opposite direction, we begin with a partition $\pi$ with largest part at most $n$ and with at most $i-1$ parts equal to $i$ for all $i\leq n$. We then arrange the parts of the partition into the shape of a DPP in order, putting as many parts in a row as possible before the part would be forced to be special, that is, before the part $i$ would be in position $(k,j)$ with $i\le j-k$. Thus $i$ will be a non-special part in the first $i$ spots in any row, that is, in positions $(k,j)$ for $k\le j< k+i$. Suppose we have filled the DPP with the parts $n, n-1, \ldots, i+1$ of $\pi$ by the above process and the result is a valid DPP with no special parts. We need to show that we can insert up to $i-1$ parts equal to $i$ and still obtain a valid DPP with no special parts.

Suppose the last part greater than $i$ was inserted in position $(k,j)$. According to our algorithm, we must insert the first $i$ in position $(k,j+1)$ if $i>j+1-k$ and in position $(k+1,k+1)$ otherwise. 
In either case, there are at least $i-2$ more columns into which we can insert $i$ as a non-special part while not violating the column inequality condition. So the columns are strictly decreasing. 

If $i$ is the smallest part in row $k$, then there will be at most $i$ parts in row $k$. This is because the $i$th part is in position $(k,k+i-1)$ and so any additional part after the $i$th part would have value less than or equal to $i$ and so would be a special part since its value would be less than or equal to $(k+i)-k$. 
Since the rows are weakly decreasing, $i$ is less than or equal to the largest part in row $k$. So we have that the number of parts of row $k$ is less than or equal to $i$ which is less than or equal to the largest part in row $k$. 
If the number of parts of row $k$ equals the largest part in row $k$, it must be that the largest part equals~$i$. So $i$ would be both the largest and smallest part in row $k$, meaning that there would be at most $i-1$ parts in row $k$ since there are at most $i-1$ parts equal to $i$ in $\pi$. Thus the number of parts in each row will always be strictly less than the largest part in that row. 

If $i$ is the largest part in row $k+1$ (where $k\geq 1$), then by our algorithm, $i$ would have been a special part if placed at the end of row $k$. That is, the first blank entry of row $k$ is in column $j$ where $i\le j-k$. So row $k$ has nonempty entries in columns $k$ through $j-1$ for a total of $(j-1)-k+1=j-k\geq i$ entries. So the number of parts in row $k$ is greater than or equal to the largest part in row $k+1$.
Therefore the number of parts in each row is greater than or equal to the largest part in the next row. 
We have verified all the conditions of a DPP, so the result is a valid DPP.


We have shown that both maps result in the desired kind of objects. The first map is the only way to make a partition from the parts of a DPP. We now show that the DPP resulting from the second map is the unique DPP with no special parts made up of exactly the parts of the partition. 
For suppose there were another DPP with no special parts whose parts made up the same partition $\pi$ with largest part at most $n$ and with at most $i-1$ parts equal to $i$ for all $i\leq n$. Let $\delta$ be the DPP produced by our algorithm and $\delta'$ be the other DPP. Look at the first row where they differ 
 and the first position of difference in that row; call it position $(k,j)$. Suppose $i$ is in position $(k,j)$ in $\delta$ and $i'$ is in position $(k,j)$ in $\delta'$ ($i'$ may be the empty part). Then either $i'<i$ or $i'$ is the empty part. 
Then there is an additional $i$ in $\delta'$ which must be put somewhere in a later row. 
Note that since $\delta$ and $\delta'$ agree up until position $(k,j)$ it must be the case that the largest part in row $k+1$ of $\delta'$ is less than or equal to $i$ (since all parts larger than $i$ occupy a position in a previous row or earlier in row $k$). So the additional $i$ in $\delta'$ must be in row $k+1$. This forces the number of elements in row $k$ of $\delta'$ to be at least $i$. So the position $(k,k+i-1)$ must not be empty in $\delta'$. In order for this entry to be a non-special part, we need its value to be greater than $(k+i-1)-k=i-1$. So the last entry in row $k$ must be greater than or equal to $i$ which forces entry $(k,j)$ to be greater than or equal to $i$. This is a contradiction. Therefore $\delta$ is the unique DPP with no special parts whose parts are those of $\pi$.
Therefore the two maps are inverses and we have a bijection.
\end{proof}

The above bijection yields an interesting generating function for DPPs with no special parts.

\begin{corollary}
The generating function for DPPs 
of order $n$ 
and no special parts (with weight equal to the sum of the parts) is 
\begin{equation}
\label{eq:permgf}
\prod_{i=1}^{n} [i]_{q^i}
\end{equation}
where $[k]_q=1+q+q^2+\ldots +q^{k-1}$ .
\end{corollary}

\begin{proof}
From basic partition theory, (\ref{eq:permgf}) is the generating function of partitions with largest part at most $n$ and at most $i-1$ parts equal to $i$ for all $i\leq n$. Since the bijection of Lemma~\ref{thm:dpppart} preserves the value of the parts, this is also the generating function for DPPs of order $n$ with no special parts. 
\end{proof}

When $q=1$ (\ref{eq:permgf}) reduces to $\prod_{i=1}^{n} i = n!$ as expected. Note that (\ref{eq:permgf}) is \emph{not} equal to $n!_q$, the typical $q$-version of $n!$. For example, for $n=3$ (\ref{eq:permgf}) becomes $(1+q^2)(1+q^3+q^6)=1+q^2+q^3+q^5+q^6+q^8$ whereas $3!_q=(1+q)(1+q+q^2)=1+2q+2q^2+q^3$.

Before proving the main theorem, we need to introduce an additional object.
\begin{definition}
Monotone triangles of order $n$
are all triangular arrays of integers with $i$~integers in row $i$ for all $1\le i\le n$, bottom row $1~2~3~\cdots~n$, 
and integer entries $a_{i,j}$ such that $a_{i,j} \le a_{i-1,j} \le a_{i,j+1} \mbox{ and } a_{i,j} < a_{i,j+1}$.
\end{definition}
It is well-known that monotone triangles of order $n$ are in bijection with $n\times n$ alternating sign matrices via the following map~\cite{BRESSOUDBOOK}. For each row of the ASM note which columns have a partial sum (from the top) of 1 in that row. Record the numbers of the columns in which this occurs 
in increasing order. This process yields a monotone triangle of order $n$.
Note that entries $a_{i,j}$ in the monotone triangle satisfying the strict diagonal inequalities \linebreak $a_{i,j} < a_{i-1,j} < a_{i,j+1}$ are in one-to-one correspondence with the $-1$ entries of the corresponding ASM.

%

We are now ready to state and prove the main theorem.

\begin{theorem}
\label{thm:bij}
There is a simple bijection between descending plane partitions
of order $n$ 
with a total of $p$ parts, $k$ parts equal to $n$, and no special parts and $n\times n$ permutation matrices with inversion number $p$ whose 1 in the last column is in row $n-k$. 
\end{theorem}

\begin{proof}
We first describe the bijection map. An example of this bijection is shown in Figure~\ref{fig:bijex}.

Begin with a DPP $\delta$
of order $n$ 
 with no special parts.
From Lemma~\ref{thm:dpppart} we know that the parts of $\delta$ form a partition with largest part at most $n$ and at most $i-1$ parts equal to $i$ for all $i\leq n$. Use these parts to make a monotone triangle of order $n$ in the following way. 
The bottom row of a monotone triangle is always $1~2~3\cdots~n$. Let $c_i$ denote the number of parts of $\delta$ equal to $i$. Beginning with $i=n$, make a continuous path (border strip) of $i$'s in the triangle starting at the $i$ in the bottom row and at each step going northeast if possible or else northwest. The path continues until there are a total of $c_i$ northwest steps in the path. In this way, the path stays as far to the east as possible and has exactly $c_i$ entries equal to their southeast diagonal neighbor.
Decrement $i$ by one and repeat until reaching $i=1$. Since there are at most $i-1$ parts equal to $i$, this process is well-defined. The resulting array is a monotone triangle of order $n$ such that there are no entries satisfying $a_{i,j} < a_{i-1,j} < a_{i,j+1}$ (i.e.\ either $a_{i,j} = a_{i-1,j}$ or $a_{i-1,j} = a_{i,j+1}$). Thus the monotone triangle
corresponds to an $n\times n$ permutation matrix $A$, since permutation matrices are alternating sign matrices with no $-1$ entries. 

\begin{figure}[htbp]
$\begin{array}{c}
\mbox{DPP}\\
\mbox{ }\\
\begin{array}{ccccc}
&&&&\\
6 & 6 & 6 & 6 & 5\\
& 5 & 4 & 4 & 4\\
&  & 3 & 3 & \\
&&&&\\
&&&&
\end{array} \end{array}
\Leftrightarrow
\begin{array}{c}
\mbox{Monotone triangle}\\
\mbox{ }\\
\begin{array}{ccccccccccc}
  & & & & & 4 & & & & & \\
  & & & & \textbf{4} & & \textbf{6} & & & & \\
  & & & 3 & & 4 & & \textbf{6} & & & \\
  & & \textbf{3} & & \textbf{4} & & \textbf{5} & & \textbf{6} & & \\
  & 1 & & \textbf{3} & & \textbf{4} & & \textbf{5} & & \textbf{6} &\\
  1 & & 2 & & 3 & & 4 & & 5 & & 6\end{array} \end{array}
\Leftrightarrow
\begin{array}{c}
\mbox{Permutation matrix}\\
\mbox{ }\\
\left( \begin{array}{cccccc}
0&0&0&1&0&0\\
0&0&0&0&0&1\\
0&0&1&0&0&0\\
0&0&0&0&1&0\\
1&0&0&0&0&0\\
0&1&0&0&0&0
\end{array}\right) \end{array}$
\caption{An example of the bijection. The bold entries in the monotone triangle are the entries equal to their southeast diagonal neighbor. These are exactly the parts of the DPP. Note that the matrix on the right represents the permutation $463512$ which has $11$ inversions. These inversions correspond to the $11$ parts of the DPP on the left.}
\label{fig:bijex}
\end{figure}

The inverse map first takes a permutation matrix $A$ to its monotone triangle.
We claim that the parts of the corresponding DPP $\delta$ are exactly the entries of the monotone triangle which are equal to their southeast diagonal neighbor, that is, entries $a_{i,j}$ such that $a_{i,j}=a_{i+1,j+1}$. Because of the shape of the monotone triangle, there are at most $i-1$ such entries equal to $i$. Thus these entries form a partition with largest entry at most $n$ and at most $i-1$ parts equal to $i$ for all $i\le n$. Using Lemma~\ref{thm:dpppart} the parts of this partition can always be formed into a unique DPP.

This is a bijection because the monotone triangle entries $a_{i,j}$ such that $a_{i,j}=a_{i+1,j+1}$ are exactly the entries coming from the northwest steps in the border strips which are exactly the entries of $\delta$.

We now show that this map takes a DPP with $p$ parts to a permutation matrix with $p$ inversions.
The inversion number of any ASM $A$ (with the matrix entry in row $i$ and column $j$ denoted $A_{ij}$) is defined as $I(A)=\sum A_{ij} A_{k\ell}$ where the sum is over all $i,j,k,\ell$ such that $i>k$ and $j<\ell$. This definition extends the usual notion of inversion in a permutation matrix. 
In~\cite{STRIKERPOSET} we noted  (using slightly different notation) 
that $I(A)$ satisfies $I(A)=E(A)+N(A)$, where $N(A)$ is the number of $-1$'s in $A$ and $E(A)$ is the number of 
entries
in the monotone triangle equal to their southeast diagonal neighbor (i.e.\ entries $a_{i,j}$ satisfying $a_{i,j}=a_{i+1,j+1}$). Since in our case, $N(A)=0$ and $E(A)$ equals the number of parts of the corresponding DPP, we have that $I(A)$ equals the number of parts of $\delta$. 

We can see that the parts of $\delta$ correspond to permutation inversions directly by noting that to convert from the monotone triangle representation of a permutation to a usual permutation $\sigma$ such that $i\rightarrow \sigma(i)$, one simply sets $\sigma(i)$ equal to the unique new entry in row $i$ of the monotone triangle. Thus for each entry of the monotone triangle $a_{i,j}$ such that $a_{i,j}=a_{i+1,j+1}$, there will be an inversion in the permutation between $a_{i,j}$ and $\sigma(i+1)$. This is because $a_{i,j}=\sigma(k)$ for some $k\leq i$ and $\sigma(k)=a_{i,j}>\sigma(i)$. These entries $a_{i,j}$ such that $a_{i,j}=a_{i+1,j+1}$ are exactly the parts of the DPP. Thus if a DPP has $p$ parts, its corresponding permutation will have $p$ inversions.

Also, observe that if the number of parts equal to $n$ in $\delta$ is $k$, then $k$ determines the position of the 1 in the last column of the permutation matrix. This is because 
the path of $n$'s in the monotone triangle must consist of exactly $k$ northwest steps (no northeast steps). So by the bijection between monotone triangles and ASMs, the $1$ in the last column of $A$ is in row $n-k$.
So finally, we have a bijection between descending plane partitions
of order $n$ 
with a total of $p$ parts, $k$ parts equal to $n$, and no special parts and $n\times n$ permutation matrices with inversion number $p$ whose 1 in the last column is in row $n-k$. 
\end{proof}
See Figure~\ref{fig:bijex} for an example of this bijection. Note that there is a direct correspondence between parts of the DPP and entries of the monotone triangle equal to their southeast neighbor. 

\section{Toward a bijection between all DPPs and ASMs}
\label{s:fullbij}

This perspective has some nice characteristics which may make the problem of finding a full bijection between ASMs and DPPs easier to approach. One such characteristic is the fact that this bijection uses monotone triangles rather than one-line permutations, which translate directly to ASMs.
Another is that there is a one-to-one correspondence between parts of the DPP and certain entries in the monotone triangle (or 
inversions of the permutation). This indicates that a full bijection should proceed by finding a one-to-one correspondence between parts of the DPP and inversions of the ASM. As discussed in the proof of Theorem~\ref{thm:bij}, the inversion number of an ASM satisfies $I(A)=E(A)+N(A)$ where $N(A)$ is the number of $-1$'s in $A$ and $E(A)$ is the number of 
entries
in the monotone triangle equal to their southeast diagonal neighbor~\cite{STRIKERPOSET}. 
Thus there should be a one-to-one correspondence between these diagonal equalities of the monotone triangle and the non-special parts of the DPP.
Difficulties arise quickly, though, since even for $n=4$ there are examples of DPPs, such as 
$\begin{array}{ccc}
4&4&3\\
&3&1
\end{array}$, whose non-special parts cannot correspond exactly to diagonal equalities of the same number in the monotone triangle. In the above example, 1 is a special part and 4 4 3 3 are non-special parts. If the parts 4 4 3 3 are each entries in the monotone triangle equal to their southeast neighbors, there is no way to fill in the rest of the entries of the monotone triangle so that there is a $-1$ in the ASM (to correspond to the special part of 1 in the DPP). Evidently, the addition of special parts to the DPP makes the relationship between non-special parts and monotone triangle diagonal equalities more complex. 

Another complicating factor is that, though there is at most one way to write any collection of numbers as a DPP with no special parts (as shown in Lemma~\ref{thm:dpppart}), the same is not true of DPPs with special parts. For example, the parts 5 5 5 3 1 can form either the DPP 
$\begin{array}{cccc}
5&5&5&1\\
&3&&
\end{array}$ or
$\begin{array}{ccc}
5&5&5\\
&3&1
\end{array}$.
Therefore, the position of the parts in the DPP matters when special parts are present (such as the 1 in this example).

Despite these difficulties, the simplicity of the bijection presented in this paper in the case of no special parts gives hope that a nice bijection exists between all DPPs and ASMs. The discovery of such a bijection would, in particular, allow us to find a weight on ASMs corresponding to the $q$-generating function (\ref{eq:qprod}) and illuminate not only the relationship between DPPs and ASMs, but also their relationships to other combinatorial objects.

\bibliographystyle{amsalpha}

\end{document}